\documentclass[12pt]{article}

\usepackage{amsmath,amssymb,amsbsy,amsfonts,amsthm,latexsym,
                     amsopn,amstext,amsxtra,euscript,amscd,stmaryrd}

\newtheorem{theorem}{Theorem}
\newtheorem{lemma}{Lemma}

\newtheorem{proposition}{Proposition}

\theoremstyle{definition}

\def\fl#1{\left\lfloor#1\right\rfloor}
\def\rf#1{\left\lceil#1\right\rceil}

\def\N{\mathbb{N}}

\def\cB{\mathcal B}

\def\cG{\mathcal G}
\def\cI{\mathcal I}
\def\cJ{\mathcal J}

\def\tO{\widetilde O}

\def\fS{\mathfrak S}

\def\Z{{\mathbb Z}}
\def\R{{\mathbb R}}

\def\eps{{\varepsilon}}
\def\e{{\rm\bf e\/}}

\def\mand{\qquad \mbox{and} \qquad}

\def\\{\cr}

\def\({\left(}

\def\){\right)}
\def\[{\left[}
\def\]{\right]}
\def\<{\langle}
\def\>{\rangle}
\def\fl#1{\left\lfloor#1\right\rfloor}
\def\rf#1{\left\lceil#1\right\rceil}
\def\le{\leqslant}
\def\ge{\geqslant}

\begin{document}

\title{\sc Representations of integers as sums of primes from a Beatty
sequence}

\author{
{\sc William D.~Banks} \\
{Department of Mathematics} \\
{University of Missouri} \\
{Columbia, MO 65211 USA} \\
{\tt bbanks@math.missouri.edu} \\
\and
{\sc Ahmet M.~G\"ulo\u glu} \\
{Department of Mathematics} \\
{University of Missouri} \\
{Columbia, MO 65211 USA} \\
{\tt ahmet@math.missouri.edu} \\
\and
{\sc C.~Wesley Nevans} \\
{Department of Mathematics} \\
{University of Missouri} \\
{Columbia, MO 65211 USA} \\
{\tt nevans@math.missouri.edu}}

\date{}

\maketitle

\newpage

\begin{abstract}
We study the problem of representing integers $N\equiv\kappa\pmod
2$ as a sum of $\kappa$ prime numbers from the Beatty sequence
$$
\cB_{\alpha,\beta}=\{n\in\N~:~n=\fl{\alpha m+\beta}\text{~for
some~}m\in\Z\},
$$
where $\alpha,\beta\in\R$ with $\alpha>1$, and $\alpha$ is
irrational and of finite type. In particular, we show that for
$\kappa=2$, almost all even numbers have such a representation if
and only if $\alpha<2$, and for any fixed integer $\kappa\ge 3$,
all sufficiently large numbers $N\equiv\kappa\pmod 2$ have such a
representation if and only if $\alpha<\kappa$.
\end{abstract}

\section{Introduction}

The celebrated 1937 theorem of Vinogradov states that every
sufficiently large odd number is the sum of three prime numbers.
However, the statement is no longer true if all three primes are
required to lie in the \emph{Beatty sequence}
$$
\cB_\alpha=\{\fl{\alpha m}~:~m\in\N\}
$$
for a fixed irrational number $\alpha>3$. Indeed, if $N$ is odd
and
\begin{equation}
\label{eq:n3} N=\fl{\alpha m_1}+\fl{\alpha m_2}+\fl{\alpha m_3}
\end{equation}
for some $m_1,m_2,m_3\in\N$, it is easy to see that
$$
N\alpha^{-1}\le m_1+m_2+m_3<N\alpha^{-1}+3\alpha^{-1}.
$$
Hence, the relation~\eqref{eq:n3} cannot hold if the fractional
part $\{N\alpha^{-1}\}$ of $N\alpha^{-1}$ lies in the open
interval $(0,1-3\alpha^{-1})$, which happens for about
$\tfrac12(1-3\alpha^{-1})X$ positive odd integers $N\le X$. On the
other hand, for an irrational number~$\alpha$ of \emph{finite
type} (see Section~\ref{sec:discr}) in the range $1<\alpha<3$, we
show that every sufficiently large odd number is the sum of three
prime numbers, each of which lies in the Beatty sequence
$\cB_\alpha$.

More generally, for fixed $\alpha,\beta\in\R$ with $\alpha>1$,
we study the problem of representing integers as sums of primes
from the \emph{non-homogeneous Beatty sequence}
$$
\cB_{\alpha,\beta}=\{n\in\N~:~n=\fl{\alpha m+\beta}\text{~for
some~}m\in\Z\}.
$$
In this paper, we prove the following:

\begin{theorem}
\label{thm:spongebob} Let $\alpha,\beta\in\R$ with $\alpha>1$, and
suppose that $\alpha$ is irrational and of finite type. Then,
\begin{itemize}
\item[$(i)$] Almost all even numbers $N$ can be expressed as the
sum of two primes from the Beatty sequence $\cB_{\alpha,\beta}$ if
and only if $\alpha<2$.

\item[$(ii)$] For any integer $\kappa\ge 3$, every sufficiently
large number $N\equiv\kappa\pmod 2$ can be expressed as the sum of
$\kappa$ primes from the Beatty sequence $\cB_{\alpha,\beta}$ if
and only if $\alpha<\kappa$.
\end{itemize}
\end{theorem}

To state our results more explicitly, we define for every integer
$\kappa\ge 2$ the function
$$
\cG_\kappa(N)=\cG_\kappa(\alpha,\beta;N)
=\sum_{\substack{n_1+\cdots+n_\kappa=N\\
n_1,\ldots,n_\kappa\in\cB_{\alpha,\beta}}}
\Lambda(n_1)\cdots\Lambda(n_\kappa)\qquad(N\ge 1),
$$
where $\Lambda$ is the \emph{von Mangoldt function}:
$$
\Lambda(n)=
\begin{cases}
\log p&\quad\text{if $n$ is a positive power of the prime $p$;} \\
0&\quad\text{otherwise.}
\end{cases}
$$
By partial summation, our estimates for $\cG_\kappa(N)$ lead to
estimates for the number of representations of an integer
$N\equiv\kappa\pmod 2$ as the sum of $\kappa$ primes from the
Beatty sequence $\cB_{\alpha,\beta}$.

Let $\psi=\psi_\alpha$ be the periodic function with period one
which is defined on the interval $(0,1]$ as follows:
\begin{equation}
\label{eq:psialphdefn} \psi(x) = \left\{  \begin{array}{ll}
1& \quad \hbox{if $0<x\le \alpha^{-1}$}; \\
0& \quad \mbox{if $\alpha^{-1}<x\le 1$}.
\end{array} \right.
\end{equation}
The function $\psi$ is closely related to the characteristic
function of the set $\cB_{\alpha,\beta}$. Let $\psi^{(1)}=\psi$,
and for every $\kappa\ge 2$, let $\psi^{(\kappa)}$ denote the
$\kappa$-fold convolution of $\psi$ with itself, defined
inductively by
$$
\psi^{(\kappa)}(x)=\int_0^1\psi^{(\kappa-1)}(x-y) \psi(y)
\,dy\qquad(\kappa\ge 2).
$$

Finally, for every $\kappa\ge 2$ we define the \emph{singular
series}
$$
\fS_\kappa(N)=\prod_{p\,\mid\,N}
\(1+\frac{(-1)^\kappa}{(p-1)^{\kappa-1}}\)
\prod_{p\,\nmid\,N}\(1+\frac{(-1)^{\kappa+1}}{(p-1)^\kappa}\)
\qquad(N\ge 1).
$$
The numbers $\fS_\kappa(N)$ arise naturally in estimates for the
number of representations of an integer as a sum of $\kappa$ prime
numbers.  Note that $\fS_\kappa(N)=0$ if and only if
$N\not\equiv\kappa\pmod 2$.

\begin{theorem}
\label{thm:main-two} Let $\alpha,\beta\in\R$ with $\alpha>1$, and
suppose that $\alpha$ is irrational and of finite type.  Then, for
any constant $C>0$, the estimate
$$
\cG_2(N)=\psi^{(2)}(\gamma N+2\delta)\fS_2(N)N+O\(\frac{N}{(\log
N)^C}\)
$$
holds for all but $O\big(X(\log X)^{-C}\big)$ integers $N\le X$,
where $\gamma=\alpha^{-1}$, $\delta=\alpha^{-1}(1-\beta)$, and the
implied constants depend only on $\alpha$ and~$C$.
\end{theorem}

\begin{theorem}
\label{thm:main-many} Let $\alpha,\beta\in\R$ with $\alpha>1$, and
suppose that $\alpha$ is irrational and of finite type.  Then, for
every integer $\kappa\ge 3$ and any constant $C>0$, the estimate
$$
\cG_\kappa(N)=\psi^{(\kappa)}(\gamma N+\kappa\delta)
\fS_\kappa(N)\frac{N^{\kappa-1}}{(\kappa-1)!}+O\(\frac{N^{\kappa-1}}{(\log
N)^C}\)
$$
holds, where $\gamma=\alpha^{-1}$, $\delta=\alpha^{-1}(1-\beta)$,
and the implied constant depends only on $\alpha$, $\kappa$ and
$C$.
\end{theorem}

The proof of Theorem~\ref{thm:main-two} is given in
Section~\ref{sec:twoprimes} (see the remark after the statement of
Theorem~\ref{thm:G2N}) and that of Theorem~\ref{thm:main-many} is
given in Section~\ref{sec:threeormoreprimes} (see the remark after
the statement of Proposition~\ref{prop:induction}).  In
Section~\ref{sec:convolutions} we study properties of the convolutions
$\psi^{(\kappa)}$ ($\kappa\ge 2$) and, in particular, derive a
sharp lower bound for values of $\psi^{(\kappa)}$ in the special
case that $\kappa=\rf{\alpha}$.  Our proof of
Theorem~\ref{thm:spongebob}, which is given in
Section~\ref{sec:proofmain}, follows immediately from the results
of Section~\ref{sec:convolutions}.

Our arguments have been strongly influenced by the treatment of
the Goldbach problem that is given in the book~\cite{IwKow} of
Iwaniec and Kowalski, and we adopt a similar notation here.  Our
underlying approach relies heavily on ideas from a recent paper of
Banks and Shparlinski~\cite{BaSh} on primes in a Beatty sequence.

\bigskip

\noindent{\bf Acknowledgements.} The authors wish to thank Ignacio
Uriarte-Tuero and Bob Vaughan for helpful conversations. This work
was done entirely at the University of Missouri-Columbia; the
support of this institution is gratefully acknowledged.

\section{Preliminaries}

\subsection{Notation}

The notation $\llbracket x\rrbracket$ is used to denote the distance
from the real number $x$ to the nearest integer; that is,
$$
\llbracket x\rrbracket=\min_{n\in\Z}|x-n|\qquad(x\in\R).
$$
We denote by $\fl{x}$, $\rf{x}$ and $\{x\}$ the greatest integer
$\le x$, the least integer $\ge x$, and the fractional part of
$x$, respectively. We also put $\e(x)=e^{2\pi ix}$ for all
$x\in\R$. As usual, we use $\mu$ and $\varphi$ to denote the
M\"obius and Euler functions, respectively.

Throughout the paper, the implied constants in symbols $O$, $\ll$
and $\gg$ may depend (where obvious) on the parameters
$\alpha,\kappa,C$ but are absolute otherwise. We recall that for
functions $F$ and $G$ the notations $F\ll G$, $G\gg F$ and
$F=O(G)$ are all equivalent to the statement that the inequality
$|F|\le c|G|$ holds for some constant $c>0$.

\subsection{Discrepancy of fractional parts}
\label{sec:discr}

Recall that the \emph{discrepancy} $D(M)$ of a sequence of (not
necessarily distinct) real numbers $a_1,a_2,\ldots,a_M\in[0,1)$ is
defined by
\begin{equation}
\label{eq:descr_defn}
D(M)=\sup_{\cI\subseteq[0,1)}\left|\frac{V(\cI,M)}{M}-|\cI|\,\right|,
\end{equation}
where the supremum is taken over all subintervals $\cI=(c,d)$ of the
interval $[0, 1)$, $V(\cI,M)$ is the number of positive integers
$m\le M$ such that $a_m\in\cI$, and $|\cI|=d-c$ is the length of
$\cI$.

For any irrational number $\gamma$ we define its \emph{type} $\tau$
by the relation
$$
\tau=\sup\Bigl\{t\in\R~:~\liminf\limits_{n\in\N}
~n^t\,\llbracket\gamma n\rrbracket=0\Bigr\}.
$$
Using~\emph{Dirichlet's approximation theorem}, it is easily seen
that $\tau\ge 1$ for every irrational number $\gamma$. The well
known theorems of Khinchin~\cite{Khin} and of Roth~\cite{Roth1,
Roth2} assert that $\tau=1$ for \emph{almost all} real numbers (in
the sense of the Lebesgue measure) and \emph{all} irrational
algebraic numbers $\gamma$, respectively; see
also~\cite{Bug,Schm}.

For every irrational number $\gamma$, it is known that the sequence
of fractional parts $\{\gamma\},\{2\gamma\},\{3\gamma\},\,\ldots\,,$
is \emph{uniformly distributed modulo~$1$} (for instance, see
\cite[Example~2.1, Chapter~1]{KuNi}). When $\gamma$ is of finite
type, this statement can be made more precise. Let
$D_{\gamma,\delta}(M)$ denote the discrepancy of the sequence of
fractional parts $(\{\gamma m + \delta\})_{m=1}^M$.
By~\cite[Theorem~3.2, Chapter~2]{KuNi} we have:

\begin{lemma}
\label{lem:discr_with_type}  Let $\gamma$ be a fixed irrational
number of finite type $\tau<\infty$.  Then, for all $\delta\in\R$
the following bound holds:
$$
D_{\gamma,\delta}(M)\le M^{-1/\tau+o(1)}\qquad(M\to\infty),
$$
where the function implied by $o(\cdot)$ depends only on $\gamma$.
\end{lemma}

\subsection{Numbers in a Beatty sequence}

The following elementary result characterizes the set of numbers
that occur in the Beatty sequence $\cB_{\alpha,\beta}$:

\begin{lemma}
\label{lem:Beatty_values} Let $\alpha,\beta\in\R$ with $\alpha>1$,
and put $\gamma=\alpha^{-1}$, $\delta=\alpha^{-1}(1-\beta)$. Then,
$n=\fl{\alpha m+\beta}$ for some integer $m$ if and only if
$0<\{\gamma n+\delta\}\le\gamma$.
\end{lemma}

\subsection{Estimates with the von Mangoldt function}

The following estimate follows immediately from the
\emph{Siegel--Walfisz theorem} (see, for example, the
book~\cite{Hux} by Huxley) using partial summation:

\begin{lemma}
\label{lem:siegel_walfisz} Let $\kappa\ge 1$ be fixed. For any
fixed constant $A>0$ and uniformly for integers $N\ge 3$ and $0\le
c<d\le(\log N)^A$ with $\gcd(c,d)=1$, the estimate
$$
\sum_{\substack{n\le N\\n\equiv c\pmod
d}}\hskip-10pt\Lambda(n)(N-n)^{\kappa-1}
=\frac{N^\kappa}{\kappa\,\varphi(d)}+O\(N^\kappa\exp\(-B(\log
N)^{1/2}\)\)
$$
holds, where $B>0$ is a constant that depends only on $\kappa$ and $A$.
\end{lemma}

We also need the following ``twisted'' version of
Lemma~\ref{lem:siegel_walfisz}:

\begin{lemma}
\label{lem:balog_perelli} Let $\kappa\ge 1$ be fixed. For an
arbitrary real number $\theta$ and coprime integers $c,d$ with
$0\le c<d$, if $|\theta-a/b|\le 1/N$ and $\gcd(a,b)=1$, then
$$
\sum_{\substack{n\le N\\n\equiv c\pmod
d}}\hskip-10pt\Lambda(n)\e(\theta
n)(N-n)^{\kappa-1}\ll\(b^{-1/2}N^\kappa
+b^{1/2}N^{\kappa-1/2}+N^{\kappa-1/5}\)(\log N)^3,
$$
where the implied constant depends only on $\kappa$.
\end{lemma}

\begin{proof}
The special case $\kappa=1$ is a simplified and weakened version
of a theorem of Balog and Perelli~\cite{BalPer} (see
also~\cite{Lav}), and the general case follows by partial
summation.
\end{proof}

\subsection{The singular series}

For every integer $\kappa\ge 2$, it is easy to check that the
singular series
$$
\fS_\kappa(N)=\prod_{p\,\mid\,N}
\(1+\frac{(-1)^\kappa}{(p-1)^{\kappa-1}}\)
\prod_{p\,\nmid\,N}\(1+\frac{(-1)^{\kappa+1}}{(p-1)^\kappa}\)
$$
satisfies the identity
\begin{equation}
\label{eq:Skappa1}
\fS_\kappa(N)=\sum_{d\,\mid\,N}\sum_{\substack{c\ge
1\\\gcd(c,d)=1}}\frac{\mu(c)^{\kappa+1}\mu(d)^\kappa
d}{\varphi(c)^\kappa\varphi(d)^\kappa},
\end{equation}
and for every $\kappa\ge 3$ we also have
\begin{equation}
\label{eq:Skappa2} \fS_\kappa(N)=\sum_{\substack{c,d\ge
1\\\gcd(d,cN)=1}}\frac{\mu(c)^\kappa\mu(d)^{\kappa+1}
d}{\varphi(c)^{\kappa-1}\varphi(d)^\kappa}.
\end{equation}
We also have the bound
\begin{equation}
\label{eq:S2Nbound} \fS_2(N)\ll\log\log N,
\end{equation}
and for every $\kappa\ge 3$,
\begin{equation}
\label{eq:SkappaNbound} \fS_\kappa(N)\ll 1.
\end{equation}

\section{Two Beatty primes}
\label{sec:twoprimes}

Fix $\alpha,\beta\in\R$ with $\alpha>1$, and suppose that $\alpha$
is irrational and of finite type. In this section, we focus our
attention on the function
$$
\cG_2(N)=\sum_{\substack{n_1+n_2=N\\n_1,n_2\in\cB_{\alpha,\beta}}}
\Lambda(n_1) \Lambda(n_2)\qquad(N\ge 1).
$$
Put $\gamma=\alpha^{-1}$ and $\delta=\alpha^{-1}(1-\beta)$, and
let $\tau$ denote the (finite) type of $\gamma$. We recall that
$\psi$ is the periodic function with period one which is defined
by~\eqref{eq:psialphdefn} on the interval $(0,1]$, and
$\psi^{(2)}=\psi*\psi$ is the convolution of $\psi$ with itself.

\begin{theorem}
\label{thm:G2N} For any complex numbers $c_N$ and any constant
$C>0$, we have
$$
\sum_{N\le X}c_N\cG_2(N)=\sum_{N\le X}c_N\psi^{(2)}(\gamma
N+2\delta) \fS_2(N)N+O\(\|c\|_2 \frac{X^{3/2}}{(\log X)^C}\),
$$
where $\|c\|_2=\(\sum_{N\le X}|c_N|^2\)^{1/2}$.
\end{theorem}

\noindent\emph{Remark.} This result immediately yields a proof of
Theorem~\ref{thm:main-two}.  Indeed, taking
$c_N=\cG_2(N)-\psi^{(2)}(\gamma N+2\delta)\fS_2(N)N$, we derive
the bound
$$
\sum_{N\le X}\big(\cG_2(N)-\psi^{(2)}(\gamma
N+2\delta)\fS_2(N)N\big)^2\ll\frac{X^3}{(\log X)^{2C}},
$$
and Theorem~\ref{thm:main-two} follows at once.

\begin{proof}[Proof of Theorem~\ref{thm:G2N}]
By Lemma~\ref{lem:Beatty_values} and the
definition~\eqref{eq:psialphdefn}, it follows that
\begin{equation}
\label{eq:G2expr} \cG_2(N)=\sum_{n_1+n_2=N}\Lambda(n_1)
\Lambda(n_2) \psi(\gamma n_1+\delta)\psi(\gamma n_2+\delta).
\end{equation}

According to a classical result of Vinogradov (see~\cite[Chapter~I,
Lemma~12]{Vin}), for any $\Delta$ such that
$$
0 < \Delta < \frac{1}{8} \mand \Delta\le
\frac{1}{2}\min\{\gamma,1-\gamma\}
$$
there is a real-valued function $\Psi$ with the following
properties:
\begin{itemize}
\item[$(i)$~~] $\Psi$ is periodic with period one;

\item[$(ii)$~~] $0 \le\Psi(x)\le 1$ for all $x\in\R$;

\item[$(iii)$~~] $\Psi(x)=\psi(x)$ if $\Delta\le \{x\}\le
\gamma-\Delta$ or if $\gamma+\Delta\le \{x\}\le 1-\Delta$;

\item[$(iv)$~~] $\Psi$ can be represented as a Fourier series:
$$
\Psi(x)=\sum_{k\in\Z}g(k)\e(kx),
$$
where $g(0)=\gamma$, and the Fourier coefficients satisfy the uniform bound
\begin{equation}
\label{eq:coeffbounds}
|g(k)|\ll\min\big\{|k|^{-1},|k|^{-2}\Delta^{-1}\big\}
\qquad(k\ne 0).
\end{equation}
\end{itemize}
{From} the properties $(i)$--$(iii)$ above, it follows that the
estimate
\begin{equation}
\label{eq:P2p2} \Psi^{(2)}(x)=\psi^{(2)}(x)+O(\Delta)
\end{equation}
holds uniformly for all $x\in\R$, where $\Psi^{(2)}$ is the
convolution $\Psi*\Psi$.

{From}~\eqref{eq:G2expr} we see that
\begin{equation}
\label{eq:basic_estimate}
\begin{split}
\cG_2(N)&=\sum_{n_1+n_2=N}\Lambda(n_1) \Lambda(n_2)
\Psi(\gamma n_1+\delta)\Psi(\gamma n_2+\delta)\\
&\qquad\qquad+O\big(V(\cI,N)(\log N)^2\big),
\end{split}
\end{equation}
where $V(\cI,N)$ is the number of positive integers $n\le N$ such
that
$$
\{\gamma n+\delta\}\in\cI=
[0,\Delta)\cup(\gamma-\Delta,\gamma+\Delta) \cup(1-\Delta,1).
$$
Since $|\cI|=4\Delta$, it follows from the
definition~\eqref{eq:descr_defn} and
Lemma~\ref{lem:discr_with_type} that
\begin{equation}
\label{eq:bound V(I,N)} V(\cI,N)\ll\Delta N+N^{1-1/(2\tau)}.
\end{equation}

Now let $K\ge\Delta^{-1}$ be a large real number (to be specified
later), and let $\Psi_K$ be the trigonometric polynomial given by
\begin{equation}
\label{eq:PKdefn} \Psi_K(x)=\sum_{|k|\le K}g(k)\e(kx).
\end{equation}
Using~\eqref{eq:coeffbounds}, we see that the estimate
\begin{equation}
\label{eq:PKP} \Psi_K(x)=\Psi(x)+O(K^{-1}\Delta^{-1})
\end{equation}
holds uniformly for all $x\in\R$, and therefore
\begin{equation}
\label{eq:PKp2est}
\Psi_K^{(2)}(x)=\Psi^{(2)}(x)+O(K^{-1}\Delta^{-1})
=\psi^{(2)}(x)+O(\Delta+K^{-1}\Delta^{-1}),
\end{equation}
where we have used~\eqref{eq:P2p2} in the second step. {From} the
definition~\eqref{eq:PKdefn} we also have
\begin{equation}
\label{eq:PK2expansion} \Psi_K^{(2)}(x)=\sum_{|k|\le
K}g(k)^2\e(kx).
\end{equation}
Inserting the estimate~\eqref{eq:PKP}
into~\eqref{eq:basic_estimate} and taking into
account~\eqref{eq:bound V(I,N)}, we derive that
\begin{equation*}
\begin{split}
\cG_2(N)&=\sum_{n_1+n_2=N}\Lambda(n_1) \Lambda(n_2)
\Psi_K(\gamma n_1+\delta)\Psi_K(\gamma n_2+\delta)\\
&\qquad\qquad+O\big(\big(\Delta+K^{-1}\Delta^{-1}+N^{-1/(2\tau)}\big)
N(\log N)^2\big).
\end{split}
\end{equation*}

For a given real number $Z\ge 2$, we now split $\Lambda(n)$ as
follows:
$$
\Lambda(n)=-\sum_{d\,\mid\,n}\mu(d)\log
d=\Lambda^\sharp(n)+\Lambda^\flat(n),
$$
where
$$
\Lambda^\sharp(n)=-\sum_{\substack{d\,\mid\,n\\d\le Z}}\mu(d)\log
d\mand\Lambda^\flat(n)=-\sum_{\substack{d\,\mid\,n\\d> Z}}\mu(d)\log
d.
$$
Then,
\begin{equation}
\label{eq:basic_estimate2}
\begin{split}
\cG_2(N)&=\cG_2^{\sharp\sharp}(N)
+2\cG_2^{\sharp\flat}(N)+\cG_2^{\flat\flat}(N)\\
&\qquad+O\big(\big(\Delta+K^{-1}\Delta^{-1}+N^{-1/(2\tau)}\big)
N(\log N)^2\big),
\end{split}
\end{equation}
where
\begin{equation*}
\begin{split}
\cG_2^{\sharp\sharp}(N)&=\sum_{n_1+n_2=N}\Lambda^\sharp(n_1)
\Lambda^\sharp(n_2)\Psi_K(\gamma n_1+\delta)\Psi_K(\gamma n_2+\delta),\\
\cG_2^{\sharp\flat}(N)&=\sum_{n_1+n_2=N}\Lambda^\sharp(n_1)
\Lambda^\flat(n_2)\Psi_K(\gamma n_1+\delta)\Psi_K(\gamma n_2+\delta),\\
\cG_2^{\flat\flat}(N)&=\sum_{n_1+n_2=N}\Lambda^\flat(n_1)
\Lambda^\flat(n_2)\Psi_K(\gamma n_1+\delta)\Psi_K(\gamma
n_2+\delta).
\end{split}
\end{equation*}

{From} now on, let $X$ be a large real parameter, and put
\begin{equation}
\label{eq:delta-K-choices1} \Delta=X^{-1/(8\tau)}\mand
K=X^{1/(4\tau)}.
\end{equation}
Then, for all $N\le X$ the estimate~\eqref{eq:basic_estimate2}
implies
$$
\cG_2(N)=\cG_2^{\sharp\sharp}(N)
+2\cG_2^{\sharp\flat}(N)+\cG_2^{\flat\flat}(N)
+O\big(X^{1-1/(10\tau)}\big).
$$
Therefore, for any complex numbers $c_N$, it follows that
\begin{equation}
\label{eq:basic_estimate3}
\begin{split}
\sum_{N\le X}c_N\cG_2(N)&=\sum_{N\le
X}c_N\big(\cG_2^{\sharp\sharp}(N)
+2\cG_2^{\sharp\flat}(N)+\cG_2^{\flat\flat}(N)\big)\\
&\qquad\qquad+O\big(\|c\|_2\,X^{3/2-1/(10\tau)}\big).
\end{split}
\end{equation}

Next, we need the following result, the proof of which is given
below:

\begin{lemma}
\label{lem:lambdaflat} For any complex numbers $u_\ell$ and $v_m$,
the bound
$$
\sum_{\ell+m+n=\fl{X}}u_\ell\,v_m\,\Lambda^\flat(n)\Psi_K(\gamma
n+\delta)\ll\|u\|_2\,\|v\|_2\,\frac{X(\log X)^2}{(\log Z)^{A}}
$$
holds with any $A>0$, where $\|u\|_2=\(\sum_{\ell\le
X}|u_\ell|^2\)^{1/2}$, $\|v\|_2=\(\sum_{m\le X}|v_m|^2\)^{1/2}$,
and the implied constant depends only on $\alpha$ and $A$.
\end{lemma}

For any complex numbers $c_N$, we have
$$
\sum_{N\le X}c_N\cG_2^{\sharp\flat}(N)=\sum_{\ell+m+n=\fl{X}}
c_{\fl{X}-\ell}\,\Lambda^\sharp(m) \Psi_K(\gamma
m+\delta)\cdot\Lambda^\flat(n)\Psi_K(\gamma n+\delta).
$$
We now apply Lemma~\ref{lem:lambdaflat} with
$$
u_\ell=\left\{
         \begin{array}{ll}
           c_{\fl{X}-\ell} & \quad\hbox{if $1\le\ell\le X$;} \\
           0 & \quad\hbox{otherwise,}
         \end{array}
       \right.
$$
and
$$
v_m=\left\{
      \begin{array}{ll}
        \Lambda^\sharp(m) \Psi_K(\gamma m+\delta)
    & \quad\hbox{if $1\le m\le X$;} \\
        0 & \quad\hbox{otherwise.}
      \end{array}
    \right.
$$
Using the trivial bound
$$
|\Lambda^\sharp(m) \Psi_K(\gamma m+\delta)|\le d(m)\log(m),
$$
where $d(m)$ is the number of positive integer divisors of $m$, it
follows that
$$
\|v\|_2^2\ll X(\log X)^5,
$$
where we have used the well known bound $\sum_{m\le X}d(m)^2\ll
X(\log X)^3$ (see, for example, the proof given by
Hua~\cite[Theorem~5.3]{Hua}; see also~\cite{Norton,Wils,Wirs}).
Hence, using Lemma~\ref{lem:lambdaflat} with $A=C+9/2$ we derive
the bound
\begin{equation}
\label{eq:sumcNGsf} \sum_{N\le X}c_N\cG_2^{\sharp\flat}(N)\ll
\|c\|_2\,\frac{X^{3/2}(\log X)^{9/2}}{(\log Z)^{C+9/2}}
\end{equation}
for any constant $C>0$. Similarly,
\begin{equation}
\label{eq:sumcNGff} \sum_{N\le
X}c_N\cG_2^{\flat\flat}(N)\ll\|c\|_2\,\frac{X^{3/2}(\log
X)^{9/2}}{(\log Z)^{C+9/2}}.
\end{equation}

Turning to the sum $\cG_2^{\sharp\sharp}(N)$,
we begin by inserting the Fourier expansion of $\Psi_K(x)$
and then changing the order of summation, obtaining
\begin{equation*}
\begin{split}
\cG_2^{\sharp\sharp}(N)&=\sum_{n_1+n_2=N}\Lambda^\sharp(n_1)
\Lambda^\sharp(n_2)
 \Psi_K(\gamma n_1+\delta)\Psi_K(\gamma n_2+\delta)\\
&=\sum_{n\le N}\Lambda^\sharp(n) \Lambda^\sharp(N-n)
 \Psi_K(\gamma n+\delta)\Psi_K(\gamma (N-n)+\delta)\\
&=\sum_{\substack{|k|\le K\\|\ell|\le K}}g(k)g(\ell)
\e(k\delta) \e(\ell(\gamma N+\delta))
\sum_{n\le N}\Lambda^\sharp(n) \Lambda^\sharp(N-n)
\e((k-\ell)\gamma n).
\end{split}
\end{equation*}
We now collect terms in double sum according to whether $k=\ell$ or not.
Writing
$$
G_2^{\sharp\sharp}(N)=\sum_{n\le N}\Lambda^\sharp(n)
\Lambda^\sharp(N-n),
$$
the contribution to $\cG_2^{\sharp\sharp}(N)$ coming
from terms with $k=\ell$ is
$$
G_2^{\sharp\sharp}(N)
\sum_{|k|\le K}g(k)^2\e(k(\gamma N+2\delta))=
\Psi_K^{(2)}(\gamma N+2\delta)G_2^{\sharp\sharp}(N),
$$
where we have used~\eqref{eq:PK2expansion} in the second step. To
bound the remainder
$$
R=\sum_{\substack{|k|,|\ell|\le K\\(k\ne\ell)}}g(k)g(\ell)
\e(k\delta) \e(\ell(\gamma N+\delta))
\sum_{n\le N}\Lambda^\sharp(n) \Lambda^\sharp(N-n)
\e((k-\ell)\gamma n),
$$
we use the following result, the proof of which is given below:

\begin{lemma}
\label{lem:lambdasharp} For every integer $k_0\ne 0$ with
$|k_0|\le 2K=2X^{1/(4\tau)}$, we have
$$
\sum_{n\le N}\Lambda^\sharp(n)\Lambda^\sharp(N-n)\e(k_0\gamma
n)\ll X^{1/2}Z^{3+4\tau},
$$
where the implied constant depends only on $\alpha$.
\end{lemma}

Using Lemma~\ref{lem:lambdasharp}, it follows that
$$
R\ll X^{1/2}Z^{3+4\tau} \sum_{|k|\le K} \big|g(k)\big|
\sum_{|\ell|\le K} \big|g(\ell)\big|\ll X^{1/2}Z^{3+4\tau}(\log
X)^2,
$$
where we have used~\eqref{eq:coeffbounds} together with our choice
of $K$.

We have therefore shown that
$$
\cG_2^{\sharp\sharp}(N)=\Psi_K^{(2)}(\gamma N+2\delta)
\,G_2^{\sharp\sharp}(N)+O\big(X^{1/2}Z^{3+4\tau} (\log X)^2\big).
$$
For any complex numbers $c_N$, it follows that
$$
\sum_{N\le X}c_N\cG_2^{\sharp\sharp}(N) =\sum_{N\le
X}c_N\Psi_K^{(2)}(\gamma N+2\delta) \,G_2^{\sharp\sharp}(N)
+O\big(\|c\|_2\,XZ^{3+4\tau}(\log X)^2\big).
$$
Now put $Z=X^{1/(9+12\tau)}$.  Using the previous estimate
together with the bounds~\eqref{eq:sumcNGsf}
and~\eqref{eq:sumcNGff}, we derive from~\eqref{eq:basic_estimate3}
the estimate
$$
\sum_{N\le X}c_N\cG_2(N)=\sum_{N\le X}c_N\Psi_K^{(2)}(\gamma
N+2\delta) \,G_2^{\sharp\sharp}(N)
+O\(\|c\|_2\,\frac{X^{3/2}}{(\log X)^C}\).
$$
Examining the proof of~\cite[Lemma~19.3]{IwKow} (which is stated
only for even numbers $N$ but holds for odd numbers as well) and
taking into account the identity~\eqref{eq:Skappa1} with
$\kappa=2$, we deduce that
$$
G_2^{\sharp\sharp}(N)=\fS_2(N)N+O\(\frac{N}{(\log N)^C}\).
$$
Using the trivial estimate
$$
\sum_{N\le X}c_N\Psi_K^{(2)}(\gamma N+2\delta)\ll\|c\|_2\,X^{1/2},
$$
it follows that
$$
\sum_{N\le X}c_N\cG_2(N)=\sum_{N\le X}c_N\Psi_K^{(2)}(\gamma
N+2\delta) \,\fS_2(N)N +O\(\|c\|_2\,\frac{X^{3/2}}{(\log X)^C}\).
$$
Finally, by~\eqref{eq:PKp2est} and our choices of $\Delta$ and
$K$, we have
$$
\Psi_K^{(2)}(x)=\psi^{(2)}(x)+O(X^{-1/(8\tau)}).
$$
In view of the trivial bound~\eqref{eq:S2Nbound}, it follows that
$$
X^{-1/(8\tau)}\sum_{N\le
X}c_N\,\fS_2(N)N\ll\|c\|_2\,X^{3/2-1/(8\tau)}\log\log X;
$$
therefore,
$$
\sum_{N\le X}c_N\cG_2(N)=\sum_{N\le X}c_N\psi^{(2)}(\gamma
N+2\delta) \,\fS_2(N)N +O\(\|c\|_2\,\frac{X^{3/2}}{(\log X)^C}\)
$$
as required.
\end{proof}

\begin{proof}[Proof of Lemma~\ref{lem:lambdaflat}]
We argue as in~\cite[Section~19.3]{IwKow} and begin with a bound
for the exponential sum
$$
S_{\Psi_K}^\flat(\xi)=\sum_{n\le X}\Lambda^\flat(n)\Psi_K(\gamma
n+\beta)\e(\xi n).
$$
{From} the definition~\eqref{eq:PKdefn}, it follows that
$$
|S_{\Psi_K}^\flat(\xi)|\le\sum_{|k|\le
K}\big|g(k)S^\flat(\xi+k\gamma)\big|,
$$
where
$$
S^\flat(\xi)=\sum_{n\le X}\Lambda^\flat(n)\e(\xi n).
$$
Using the bound~(19.17) from~\cite{IwKow} together
with~\eqref{eq:coeffbounds}, we immediately deduce that the
uniform bound
\begin{equation}
\label{eq:Spkfxibound} |S_{\Psi_K}^\flat(\xi)|\ll \frac{X\log X\log
K}{(\log Z)^A}\qquad(\xi\in\R)
\end{equation}
holds with any fixed constant $A>0$.

To complete the proof, we observe that
\begin{equation*}
\begin{split}
&\sum_{\ell+m+n=\fl{X}}u_\ell\,v_m\,
\Lambda^\flat(n)\Psi_K(\gamma n+\delta)\\
&\qquad\qquad=\int_0^1\Big(\sum_{\ell\le
X}u_\ell\,\e(\xi\ell)\Big)\Big(\sum_{m\le X}v_m\,\e(\xi
m)\Big)S_{\Psi_K}^\flat(\xi)\e(-\xi\fl{X})\,d\xi.
\end{split}
\end{equation*}
Applying the Cauchy-Schwarz inequality and
using~\eqref{eq:Spkfxibound} (with $K=X^{1/(4\tau)}$) together
with the equalities
$$
\int_0^1\left|\sum_{\ell\le X}u_\ell\,\e(\xi\ell)\right|^2\,d\xi
=\sum_{\ell\le X}|u_\ell|^2
$$
and
$$
\int_0^1\left|\sum_{m\le X}v_m\,\e(\xi m)\right|^2\,d\xi =\sum_{m\le
X}|v_m|^2,
$$
we obtain the stated bound.
\end{proof}

\begin{proof}[Proof of Lemma~\ref{lem:lambdasharp}]
We have:
\begin{equation}
\label{eq:spiderman}
\begin{split}
&\sum_{n\le N}\Lambda^\sharp(n)\Lambda^\sharp(N-n)\e(k_0\gamma
n)\\
&\qquad=\sum_{n\le N}\biggl(\,\sum_{\substack{d_1\,\mid\,n\\d_1\le
Z}}\mu(d_1)\log
d_1\biggl)\biggl(\,\sum_{\substack{d_2\,\mid\,N-n\\d_2\le
Z}}\mu(d_2)\log d_2\biggl)\e(k_0\gamma n)\\
&\qquad=\sum_{d_1,d_2\le Z}\mu(d_1)\mu(d_2)\log d_1\log
d_2\sum_{\substack{\ell_1,\ell_2\ge
1\\\ell_1d_1+\ell_2d_2=N}}\e(k_0\gamma\ell_1d_1).
\end{split}
\end{equation}
If $\ell_1\ge 1$, then $\ell_1d_1+\ell_2d_2=N$ for some $\ell_2\ge
1$ if and only if $\ell_1<N/d_1$, $f=\gcd(d_1,d_2)$ is a divisor of
$N$, and
$$
\ell_1(d_1/f)\equiv (N/f)\pmod{d_2/f}.
$$
Let $a$ be the least positive integer such that
$$
a\equiv(d_1/f)^{-1}(N/f)\pmod{d_2/f}
$$
Therefore, $\ell_1$ varies over the set
$\{a,a+d_2/f,\ldots,a+(L-1)d_2/f\}$, where
$$
L=\fl{\frac{N/d_1-a}{d_2/f}}=\frac{N}{[d_1,d_2]}+O(1),
$$
and it follows that
\begin{equation}
\begin{split}
\label{eq:superman}
\sum_{\substack{\ell_1,\ell_2\ge
1\\\ell_1d_1+\ell_2d_2=N}}\e(k_0\gamma\ell_1d_1)&=\e(k_0\gamma
ad_1)\sum_{j=0}^{L-1}\e(k_0\gamma j\,[d_1,d_2])\\
&\ll\frac{1}{\big\llbracket k_0\gamma\, [d_1,d_2]\big\rrbracket},
\end{split}
\end{equation}
where we have used a standard estimate in the second step (see,
for example, \cite[Chapter~1, Lemma~1]{Koro}). Since $\gamma$ is
of type $\tau$, we have
$$
\llbracket\gamma n\rrbracket\gg n^{-2\tau}\qquad(n\ge 1),
$$
where the implied constant depends on $\alpha$; thus,
$$
\frac{1}{\big\llbracket k_0\gamma\, [d_1,d_2]\big\rrbracket}\ll
k_0^{2\tau}\, [d_1,d_2]^{2\tau} \le
(2X^{1/(4\tau)})^{2\tau}Z^{4\tau} \ll X^{1/2}Z^{4\tau}.
$$
Combining this bound with~\eqref{eq:spiderman} and~\eqref{eq:superman},
and using the trivial bound
$$
\sum_{d_1,d_2\le Z}\log d_1\log d_2\le Z^2(\log Z)^2\ll Z^3,
$$
we obtain the desired result.
\end{proof}

\section{Three or more Beatty primes}
\label{sec:threeormoreprimes}

In what follows, we use the same notation as in the proof of
Theorem~\ref{thm:G2N}, except that we now define
$$
\Delta=N^{-1/(8\tau)}\mand K=N^{1/(4\tau)}
$$
instead of~\eqref{eq:delta-K-choices1}.  With these choices, we
have the following analog of~\eqref{eq:PKp2est} for every
$\kappa\ge 2$:
\begin{equation}
\label{eq:PKpkest} \Psi_K^{(\kappa)}(x)
=\psi^{(\kappa)}(x)+O(N^{-1/(8\tau)})\qquad(x\in\R).
\end{equation}
Also,
\begin{equation}
\label{eq:PKkexpansion} \Psi_K^{(\kappa)}(x)=\sum_{|\ell|\le
K}g(\ell)^\kappa\e(\ell x).
\end{equation}

\begin{proposition}
\label{prop:induction} Let $\kappa\ge 2$ be fixed. If, for any
constant $C>0$, the estimate
\begin{equation}
\label{eq:curious1} \cG_\kappa(n)=\Psi_K^{(\kappa)}(\gamma
n+\kappa\delta)
\fS_\kappa(n)\frac{n^{\kappa-1}}{(\kappa-1)!}+O\(\frac{n^{\kappa-1}}{(\log
n)^C}\)
\end{equation}
holds for all but $O\big(N(\log N)^{-C}\big)$ integers $n\le N$,
then the estimate
\begin{equation}
\label{eq:curious2} \cG_{\kappa+1}(N)=\Psi_K^{(\kappa+1)}(\gamma
N+(\kappa+1)\delta)
\fS_{\kappa+1}(N)\frac{N^\kappa}{\kappa!}+O\(\frac{N^\kappa}{(\log
N)^C}\)
\end{equation}
holds with any constant $C>0$.
\end{proposition}

\noindent\emph{Remark.} This result immediately yields a proof of
Theorem~\ref{thm:main-many}.  Indeed, using~\eqref{eq:S2Nbound}
and~\eqref{eq:PKpkest} we obtain~\eqref{eq:curious1} with
$\kappa=2$.  By induction, Proposition~\ref{prop:induction}
implies that~\eqref{eq:curious2} holds for every fixed $\kappa\ge
2$. Replacing $\kappa$ by $\kappa-1$ in~\eqref{eq:curious2} and
then using the estimate~\eqref{eq:PKpkest} again, we obtain the
statement of Theorem~\ref{thm:main-many}.

\begin{proof}[Proof of Proposition~\ref{prop:induction}]
To simplify our exposition in what follows, for any functions
$F=F(N)$ and $G=G(N)$ we use notation
$$
F=\tO(G)
$$
to mean that for any choice of the constant $C>0$ the inequality
$$
|F|\le c\,\frac{|G|}{(\log N)^C}
$$
holds for all $N\ge 2$ with a constant $c>0$ that depends only on
$\alpha$, $\kappa$ and $C$.

By Lemma~\ref{lem:Beatty_values} and the
definition~\eqref{eq:psialphdefn}, we have
\begin{equation*}
\begin{split}
\cG_{\kappa+1}(N)&=\sum_{n_1+\cdots+n_{\kappa+1}=N}
\Lambda(n_1)\cdots\Lambda(n_{\kappa+1}) \psi(\gamma
n_1+\delta)\cdots\psi(\gamma n_{\kappa+1}+\delta)\\
&=\sum_{n\le N}\Lambda(N-n)\psi(\gamma(N-n)+\delta)\cG_\kappa(n)\\
&={\sum_{n\le
N}}^*\Lambda(N-n)\psi(\gamma(N-n)+\delta)\cG_\kappa(n)
+\tO(N^\kappa),
\end{split}
\end{equation*}
where ${\sum}^*$ indicates that the sum is restricted to integers
$n$ satisfying~\eqref{eq:curious1}; note that we have used the
trivial bound
$$
\Lambda(N-n)\psi(\gamma(N-n)+\delta)\cG_\kappa(n)\ll
N^{\kappa-1}(\log N)^\kappa
$$
to estimate the contribution from exceptional integers.
By~\eqref{eq:curious1}, the previous sum is equal to
$$
{\sum_{n\le
N}}^*\Lambda(N-n)\psi(\gamma(N-n)+\delta)\Psi_K^{(\kappa)}(\gamma
n+\kappa\delta) \fS_\kappa(n)\frac{n^{\kappa-1}}{(\kappa-1)!}
+\tO(N^\kappa).
$$
We now extend the sum to all integers $n\le N$,
using~\eqref{eq:S2Nbound} or~\eqref{eq:SkappaNbound} to bound
$\fS_\kappa(n)$ for each exceptional $n$, then we replace $\psi$
with $\Psi_K$ using~\eqref{eq:PKpkest} to control the error term.
Finally, replacing $n$ by $N-n$, we see that $\cG_{\kappa+1}(N)$
is equal to
$$
\sum_{n\le N}\Lambda(n)\Psi_K(\gamma
n+\delta)\Psi_K^{(\kappa)}(\gamma(N-n)+\kappa\delta)
\fS_\kappa(N-n)\frac{(N-n)^{\kappa-1}}{(\kappa-1)!}
+\tO(N^\kappa).
$$
In this sum, we substitute the Fourier
expansions~\eqref{eq:PKdefn} and~\eqref{eq:PKkexpansion} for
$\Psi_K$ and $\Psi_K^{(\kappa)}$, respectively, then change the
order of summation, obtaining
\begin{equation}
\label{eq:G-etc} \cG_{\kappa+1}(N)=\sum_{\substack{|k|\le
K\\|\ell|\le K}}g(k)g(\ell)^\kappa\e(k\delta+\ell\gamma
N+\ell\kappa\delta)\frac{S_{k,\ell}(N)}{(\kappa-1)!}
+\tO(N^\kappa),
\end{equation}
where
$$
S_{k,\ell}(N)=\sum_{n\le N}\Lambda(n)\e((k-\ell)\gamma n)
\fS_\kappa(N-n)(N-n)^{\kappa-1}.
$$

We now show that the main contribution to $\cG_{\kappa+1}(N)$
comes from the sums $S_{k,\ell}(N)$ with $k=\ell$.  To this end,
we use~\eqref{eq:Skappa1} to write
\begin{equation*}
\begin{split}
S_{k,\ell}(N)&=\sum_{n\le N}\Lambda(n)\e((k-\ell)\gamma n)
(N-n)^{\kappa-1}\sum_{d\,\mid\,N-n}\sum_{\substack{c\ge
1\\\gcd(c,d)=1}}\frac{\mu(c)^{\kappa+1}\mu(d)^\kappa
d}{\varphi(c)^\kappa\varphi(d)^\kappa}\\
&=\sum_{d\le N}\sum_{\substack{c\ge
1\\\gcd(c,d)=1}}\frac{\mu(c)^{\kappa+1}\mu(d)^\kappa
d}{\varphi(c)^\kappa\varphi(d)^\kappa}\,T_{k,\ell,d}(N),
\end{split}
\end{equation*}
where
$$
T_{k,\ell,d}(N)=\sum_{\substack{n\le N\\n\equiv N\pmod
d}}\Lambda(n)\e((k-\ell)\gamma n)(N-n)^{\kappa-1}.
$$
Using the trivial uniform bound
$$
T_{k,\ell,d}(N)\ll\frac{N^\kappa\log N}{d}
$$
and the well known lower bound $\varphi(d)\gg d/\log\log d$, we
have for any $y>3$ (since $\kappa\ge 2$):
\begin{equation*}
\begin{split}
\sum_{d>y}\sum_{\substack{c\ge
1\\\gcd(c,d)=1}}\frac{\mu(c)^{\kappa+1}\mu(d)^\kappa
d}{\varphi(c)^\kappa\varphi(d)^\kappa}\,T_{k,\ell,d}(N)
&\ll\sum_{d>y}\frac{d(\log\log
d)^\kappa}{d^\kappa}\,\frac{N^\kappa\log N}{d}\\
&\ll N^\kappa\log N\sum_{d>y}\frac{1}{d^{3/2}}
\ll\frac{N^\kappa\log N}{y^{1/2}}.
\end{split}
\end{equation*}
Taking $y=(\log N)^A$ with $A=2C+2$ and $C>0$ arbitrary, we derive
that
$$
S_{k,\ell}(N)=\sum_{d\le (\log N)^A}\sum_{\substack{c\ge
1\\\gcd(c,d)=1}}\frac{\mu(c)^{\kappa+1}\mu(d)^\kappa
d}{\varphi(c)^\kappa\varphi(d)^\kappa}\,T_{k,\ell,d}(N)
+O\(\frac{N^\kappa}{(\log N)^C}\).
$$
Next, we observe that if $d\le(\log N)^A$ and $\gcd(d,N)\ne 1$,
then the number $\omega(d)$ of distinct prime divisors of $d$
satisfies the bound $\omega(d)\ll\log\log N$, and it is easy to
see that the bound
$$
T_{k,\ell,d}(N)\ll N^{\kappa-1}\log N\log\log N
$$
holds for all such $d$. Using this estimate in the preceding
expression for $S_{k,\ell}(N)$, it follows that
$$
S_{k,\ell}(N)=\sum_{\substack{d\le (\log
N)^A\\\gcd(d,N)=1}}\sum_{\substack{c\ge
1\\\gcd(c,d)=1}}\frac{\mu(c)^{\kappa+1}\mu(d)^\kappa
d}{\varphi(c)^\kappa\varphi(d)^\kappa}\,T_{k,\ell,d}(N)
+O\(\frac{N^\kappa}{(\log N)^C}\).
$$

In the case that $k=\ell$, Lemma~\ref{lem:siegel_walfisz}
immediately implies that
$$
T_{k,k,d}(N)=\sum_{\substack{n\le N\\n\equiv N\pmod
d}}\Lambda(n)(N-n)^{\kappa-1}
=\frac{N^\kappa}{\kappa\,\varphi(d)}+\tO(N^\kappa),
$$
and therefore,
$$
S_{k,k}(N)=\frac{N^\kappa}{\kappa}\sum_{\substack{d\le (\log
N)^A\\\gcd(d,N)=1}}\sum_{\substack{c\ge
1\\\gcd(c,d)=1}}\frac{\mu(c)^{\kappa+1}\mu(d)^\kappa
d}{\varphi(c)^\kappa\varphi(d)^{\kappa+1}}
+O\(\frac{N^\kappa}{(\log N)^C}\).
$$
Since
$$
\sum_{\substack{d>(\log N)^A\\\gcd(d,N)=1}}\sum_{\substack{c\ge
1\\\gcd(c,d)=1}}\frac{\mu(c)^{\kappa+1}\mu(d)^\kappa
d}{\varphi(c)^\kappa\varphi(d)^{\kappa+1}}\ll\sum_{d>(\log
N)^A}\frac{(\log\log d)^{\kappa+1}}{d^\kappa}\ll\frac{1}{(\log
N)^{C+1}},
$$
and $C$ is arbitrary, it follows that
$$
S_{k,k}(N)=\frac{N^\kappa}{\kappa}\sum_{\substack{c,d\ge
1\\\gcd(d,cN)=1}}\frac{\mu(c)^{\kappa+1}\mu(d)^\kappa
d}{\varphi(c)^\kappa\varphi(d)^{\kappa+1}}+\tO(N^\kappa).
$$
Finally, using~\eqref{eq:Skappa2} (with $\kappa$ replaced by
$\kappa+1$) we deduce that
\begin{equation}
\label{eq:SkkNestx}
S_{k,k}(N)=\fS_{\kappa+1}(N)\frac{N^\kappa}{\kappa}
+\tO(N^\kappa)\qquad(|k|\le K).
\end{equation}

To treat the case $k\ne\ell$, we use the following result, the
proof of which is given below:

\begin{lemma}
\label{lem:needsname} There exists a constant $\eta>0$ that
depends only on $\alpha$ with the following property.  For any
positive integer $d$ coprime to $N$, and any nonzero integer $k_0$
such that $|k_0|\le 2N^{1/(4\tau)}$, the bound
$$
\sum_{\substack{n\le N\\n\equiv N\pmod d}}\Lambda(n)\e(k_0\gamma
n)(N-n)^{\kappa-1}\ll N^{\kappa-\eta}
$$
holds, where the implied constant depends only on $\kappa$.
\end{lemma}

By Lemma~\ref{lem:needsname} we have for all $|k|,|\ell|\le K$
with $k\ne\ell$:
$$
T_{k,\ell,d}(N)=\sum_{\substack{n\le N\\n\equiv N\pmod
d}}\Lambda(n)\e((k-\ell)\gamma n)(N-n)^{\kappa-1}=\tO(N^\kappa),
$$
and therefore,
\begin{equation}
\label{eq:SklNestx}
S_{k,\ell}(N)=\tO(N^\kappa)\qquad(|k|,|\ell|\le K,~k\ne\ell).
\end{equation}

Inserting the estimates~\eqref{eq:SkkNestx}
and~\eqref{eq:SklNestx} into~\eqref{eq:G-etc}, and taking into
account~\eqref{eq:coeffbounds}, it follows that
\begin{equation*}
\begin{split} \cG_{\kappa+1}(N)&=\fS_{\kappa+1}(N)
\frac{N^\kappa}{\kappa!}\sum_{|k|\le
K}g(k)^{\kappa+1}\e(k\gamma N+(\kappa+1)k\delta)+\tO(N^\kappa)\\
&=\Psi_K^{(\kappa+1)}(\gamma N+(\kappa+1)\delta)
\fS_{\kappa+1}(N)\frac{N^\kappa}{\kappa!}+\tO(N^\kappa),
\end{split}
\end{equation*}
and this completes the proof.
\end{proof}

\begin{proof}[Proof of Lemma~\ref{lem:needsname}]
Fix a constant $\varrho$ such that
$$
1\le \tau<\varrho<2\tau.
$$
Since $\gamma$ is of type $\tau$, for some constant $c_0>0$ we
have
\begin{equation}
\label{eq:Liou} \llbracket\gamma m\rrbracket>c_0
m^{-\varrho}\qquad(m\ge 1).
\end{equation}
Taking $c_0$ smaller if necessary, we can assume that
$c_0<2^\varrho$. Put
$$
c_1=2^\varrho/c_0\mand\eps=1/(4\tau+2).
$$

Let $d$ and $k_0$ be integers with the properties stated in the
lemma; without loss of generality, we can assume that $k_0$ is
positive. Let $a/b$ be the convergent in the continued fraction
expansion of $k_0\gamma$ that has the largest denominator $b$ not
exceeding $c_1N^{1-\eps}$; then,
\begin{equation}
\label{eq:theta-mix1} \left|k_0\gamma-\frac{a}{b}\right|
\le\frac{1}{bc_1 N^{1-\eps}}=\frac{c_0}{b2^\varrho N^{1-\eps}}.
\end{equation}
Multiplying by $b$ and taking~\eqref{eq:Liou} into account, we
have
$$
\frac{c_0}{2^\varrho N^{1-\eps}}\ge\left|bk_0\gamma-a\right|
\ge\llbracket bk_0\gamma\rrbracket>c_0 (bk_0)^{-\varrho}.
$$
Thus, since $k_0\le 2N^{1/(4\tau)}$ and $\varrho<2\tau$, it
follows that
\begin{equation}
\label{eq:theta-mix2} b\ge
N^{(1-\eps)/(2\tau)-1/(4\tau)}=N^{\eps}.
\end{equation}
Inserting~\eqref{eq:theta-mix2} into~\eqref{eq:theta-mix1} and
recalling that $c_0<2^\varrho$, we conclude  that
$$
\left|k_0\gamma-\frac{a}{b}\right|\le\frac{1}{N}.
$$
We are therefore in a position to apply
Lemma~\ref{lem:balog_perelli} with $\theta=k_0\gamma$, and this
yields the stated result immediately since $N^\eps\le b\le
c_1N^{1-\eps}$.
\end{proof}

\section{Convolutions with $\psi$}
\label{sec:convolutions}

In this section, we focus on properties of the $\kappa$-fold
convolutions of $\psi$.  We recall that $\psi$ is the periodic
function with period one defined by
$$
\psi(x) = \left\{
\begin{array}{ll}
1& \quad \hbox{if $0<\{x\}\le \gamma$}; \\
0& \quad \mbox{if $\gamma<\{x\}<1$ or $x\in\Z$}.
\end{array} \right.
$$
We assume that $\gamma=\alpha^{-1}<1$. As before, we put
$\psi^{(1)}=\psi$, and for every $\kappa\ge 2$, we denote by
$\psi^{(\kappa)}$ the $\kappa$-fold convolution of $\psi$ with
itself:
$$
\psi^{(\kappa)}(x)=\int_0^1\psi^{(\kappa-1)}(x-y)\psi(y)\,dy
=\int_{x-\gamma}^x\psi^{(\kappa-1)}(y)\,dy.
$$
Since $0\le\psi(x)\le\gamma$ for all $x\in\R$, it is easy to see
that
$$
0\le\psi^{(\kappa)}(x)\le\gamma^{\kappa-1}\qquad(\kappa\ge
1,~x\in\R).
$$
Note that $\psi^{(\kappa)}$ is continuous for $\kappa\ge 2$ and
differentiable for $\kappa\ge 3$.

\begin{proposition}
\label{prop:ducktape} If $\kappa\ge\rf{\alpha}$, then there exists
a constant $c>0$ which depends only on $\alpha$ and $\kappa$ such
that $\psi^{(\kappa)}(x)\ge c$ for all $x\in\R$.
\end{proposition}

\begin{proof}
By periodicity, it suffices to prove this for all $x$ in
$[\eps,1+\eps]$ for some $\eps>0$.  Since
$\kappa\gamma\ge\rf{\alpha}/\alpha>1$, there exists $\eps>0$ such
that $1+2\eps\le\kappa\gamma$. Fixing $\eps$, it is easy to see
that for every $x\in[\eps,1+\eps]$ the closed intervals
$$
\cI_x=\left[\frac{x}{\kappa}-\frac{\eps}{\kappa}\,,
\frac{x}{\kappa}+\frac{\eps}{\kappa}\right]\mand
\cJ_x=\left[\frac{x}{\kappa}-\frac{\eps}{\kappa(\kappa-1)}\,,
\frac{x}{\kappa}+\frac{\eps}{\kappa(\kappa-1)}\right]
$$
are contained in $[0,\gamma]$. Also, if $y_j\in\cJ_x$ for
$j=1,\ldots,\kappa-1$, then the number $x-y_1-\cdots-y_{\kappa-1}$
lies in $\cI_x$. Therefore,
\begin{equation*}
\begin{split}
\psi^{(\kappa)}(x)&=\int_0^1\cdots\int_0^1\psi(y_1)\cdots\psi(y_{\kappa-1})
\psi(x-y_1-\cdots-y_{\kappa-1})\,dy_1\cdots dy_{\kappa-1}\\
&\ge\int\limits_{\cJ_x}\cdots\int\limits_{\cJ_x}\,dy_1\cdots
dy_{\kappa-1}=\(\frac{2\eps}{\kappa(\kappa-1)}\)^{\kappa-1}
\end{split}
\end{equation*}
for all $x\in[\eps,1+\eps]$.
\end{proof}

The remainder of this section is devoted to the problem of finding
a sharp lower bound for $\psi^{(\kappa)}(x)$ in the special case
that $\kappa=\rf{\alpha}$, which is given in
Theorem~\ref{thm:sharpbound} below.

\begin{lemma}
\label{lem:symmetry} If $\kappa\ge 2$, then
$\psi^{(\kappa)}(x)=\psi^{(\kappa)}(\kappa\gamma-x)$ for all
$x\in\R$.
\end{lemma}

\begin{proof}
Let $\psi_0$ be the characteristic function of the set of real
numbers $x$ such that $\llbracket x\rrbracket\le\gamma/2$.
Clearly, $\psi(x)=\psi_0(x-\gamma/2)$ for all $x\in\R\setminus\Z$,
and by induction on $\kappa$, we have
$\psi^{(\kappa)}(x)=\psi_0^{(\kappa)}(x-\kappa\gamma/2)$ for all
$\kappa\ge 2$ and $x\in\R$.  Since $\psi_0$ is an even function,
so is $\psi_0^{(\kappa)}$ for all $\kappa\ge 2$; therefore,
$$
\psi^{(\kappa)}(x)=\psi_0^{(\kappa)}(x-\kappa\gamma/2)
=\psi_0^{(\kappa)}(\kappa\gamma/2-x)=\psi^{(\kappa)}(\kappa\gamma-x)
$$
for all $\kappa\ge 2$ and $x\in\R$.
\end{proof}

\begin{lemma}
\label{lem:vanishing} If $1\le\kappa<\rf{\alpha}$ and
$x\in(\kappa\gamma,1]$, then $\psi^{(\kappa)}(x)=0$.
\end{lemma}

\begin{proof}
When $\kappa=1$, this follows from the definition of $\psi$. Now
suppose that $\psi^{(\kappa-1)}(x)=0$ for all
$x\in((\kappa-1)\gamma,1]$, where $\kappa\ge 2$.  Then, for each
$x\in(\kappa\gamma,1]$ the interval $[x-\gamma,x]$ is contained in
$((\kappa-1)\gamma,1]$; therefore,
$$
\psi^{(\kappa)}(x)=\int_{x-\gamma}^x\psi^{(\kappa-1)}(y)\,dy=0,
$$
and the result follows by induction.
\end{proof}

The next result is an easy consequence of
Lemma~\ref{lem:vanishing}:

\begin{lemma}
\label{lem:newintegral} If $2\le\kappa<\rf{\alpha}$ and
$x\in[0,\gamma]$, then
$$
\psi^{(\kappa)}(x)=\int_0^x\psi^{(\kappa-1)}(y)\,dy.
$$
The same result holds for $\kappa=\rf{\alpha}$ and
$x\in[\kappa\gamma-1,\gamma]$.
\end{lemma}

\begin{lemma}
\label{lem:exact} For $1 \le \kappa<\rf{\alpha}$ and
$x\in(0,\gamma]$, we have
$$
\psi^{(\kappa)}(x)=\frac{x^{\kappa-1}}{(\kappa-1)!}.
$$
\end{lemma}

\begin{proof}
This is immediate for $\kappa=1$. Suppose that
$\psi^{(\kappa-1)}(x)=x^{\kappa-2}/(\kappa-2)!$ for
$x\in(0,\gamma]$, where $2\le\kappa<\rf{\alpha}$.  Then, by
Lemma~\ref{lem:newintegral} we have
$$
\psi^{(\kappa)}(x)=\int_0^x\psi^{(\kappa-1)}(y)\,dy=\int_0^x
\frac{y^{\kappa-2}}{(\kappa-2)!}\,dy=\frac{x^{\kappa-1}}{(\kappa-1)!},
$$
and the result follows by induction.
\end{proof}

\begin{lemma}
\label{lem:increasing} If $1\le\kappa<\rf{\alpha}$, then
$\psi^{(\kappa)}$ is increasing on $[0,\kappa\gamma/2]$.
\end{lemma}

\begin{proof}
For $\kappa=1$ this is immediate, and for $\kappa=2$, it follows
from the fact that $\psi^{(2)}(x)=x$ for $x\in[0,\gamma]$ by
Lemma~\ref{lem:exact} and the continuity of $\psi^{(2)}$. Now
suppose that $\psi^{(\kappa-1)}$ is increasing on
$[0,(\kappa-1)\gamma/2]$, where $\kappa\ge 3$. Since
$\psi^{(\kappa)}$ is differentiable, we have for
$x\in[\gamma,(\kappa-1)\gamma/2]$:
$$
\frac{d\psi^{(\kappa)}(t)}{dt}\biggl|_{t=x}
=\psi^{(\kappa-1)}(x)-\psi^{(\kappa-1)}(x-\gamma)\ge 0.
$$
If $x\in[0,\gamma]$, then by Lemma~\ref{lem:newintegral} it
follows that
$$
\frac{d\psi^{(\kappa)}(t)}{dt}\biggl|_{t=x}
=\psi^{(\kappa-1)}(x)-\psi^{(\kappa-1)}(0)\ge 0.
$$
Finally, suppose that $x\in[(\kappa-1)\gamma/2,\kappa\gamma/2]$.
Since $\psi^{(\kappa-1)}$ is increasing on
$[0,(\kappa-1)\gamma/2]$, it is decreasing on
$[(\kappa-1)\gamma/2,(\kappa-1)\gamma]$ by
Lemma~\ref{lem:symmetry}; therefore, using the same lemma we have
\begin{equation*}
\begin{split}
\frac{d\psi^{(\kappa)}(t)}{dt}\biggl|_{t=x}
&=\psi^{(\kappa-1)}(x)-\psi^{(\kappa-1)}(x-\gamma)\\
&\ge\psi^{(\kappa-1)}(\kappa\gamma/2)-
\psi^{(\kappa-1)}((\kappa-2)\gamma/2)=0,
\end{split}
\end{equation*}
and the proof is completed by induction.
\end{proof}

\begin{theorem}
\label{thm:sharpbound} For $\kappa=\rf{\alpha}$, the sharp lower
bound
\begin{equation*}
\psi^{(\kappa)} (x) \ge \frac{(\kappa\gamma
-1)^{\kappa-1}}{2^{\kappa-2}(\kappa-1)!}
\end{equation*}
holds uniformly for all $x\in\R$.
\end{theorem}

\begin{proof}
Since $\psi^{(\kappa)}$ has period one, we can assume that $x\in
[0,1]$.

Using Lemmas~\ref{lem:symmetry} and~\ref{lem:newintegral} and
arguing as in the proof of Lemma~\ref{lem:increasing}, one sees
that $\psi^{(\kappa)}$ is increasing on the interval
$[\kappa\gamma-1,\kappa\gamma/2]$ and decreasing on the interval
$[\kappa\gamma/2,1]$.  Therefore,
$$
\psi^{(\kappa)}(x)\ge\psi^{(\kappa)}(1)=\psi^{(\kappa)}(0)
$$
for all $x\in[\kappa\gamma-1,1]$.  On the other hand, for
$x\in[0,\kappa\gamma-1]$ we have by Lemmas~\ref{lem:symmetry},
\ref{lem:vanishing} and~\ref{lem:exact}:
\begin{equation*}
\begin{split}
\psi^{(\kappa)}(x)&=\int_{x+1-\gamma}^1\psi^{(\kappa-1)}(y)\,dy
+\int_0^x\psi^{(\kappa-1)}(y)\,dy\\
&=\int_{(\kappa-1)\gamma-1}^{\kappa\gamma-1-x}\psi^{(\kappa-1)}(y)\,dy
+\int_0^x\psi^{(\kappa-1)}(y)\,dy\\
&=\int_0^{\kappa\gamma-1-x}\psi^{(\kappa-1)}(y)\,dy
+\int_0^x\psi^{(\kappa-1)}(y)\,dy=f(x),
\end{split}
\end{equation*}
where
$$
f(x)=\frac{(\kappa\gamma-1-x)^{\kappa-1}
+x^{\kappa-1}}{(\kappa-1)!}.
$$
Since the function $f(x)$ attains its minimum on
$[0,\kappa\gamma-1]$ at $x=(\kappa\gamma-1)/2$, we obtain the
stated result.
\end{proof}

\section{Proof of Theorem~\ref{thm:spongebob}}
\label{sec:proofmain}

Suppose that $\kappa<\alpha$.  If $N\equiv\kappa\pmod 2$, and
\begin{equation}
\label{eq:nkappa} N=\fl{\alpha m_1+\beta}+\fl{\alpha
m_2+\beta}+\cdots+\fl{\alpha m_\kappa+\beta}
\end{equation}
for some $m_1,\ldots,m_\kappa\in\N$, then
$$
(N-\kappa\beta)\alpha^{-1}\le
m_1+\cdots+m_\kappa<(N-\kappa\beta)\alpha^{-1}+\kappa\alpha^{-1}.
$$
Therefore, the relation~\eqref{eq:nkappa} cannot hold if the
fractional part $\{(N-\kappa\beta)\alpha^{-1}\}$ of
$(N-\kappa\beta)\alpha^{-1}$ lies in the open interval
$(0,1-\kappa\alpha^{-1})$, which happens for about
$\tfrac12(1-\kappa\alpha^{-1})X$ positive integers $N\le X$ with
$N\equiv\kappa\pmod 2$.  This proves the forward implications of
the statements in Theorem~\ref{thm:spongebob}. The reverse
implications follow immediately from Theorems~\ref{thm:main-two}
and~\ref{thm:main-many} combined with the lower bound of
Proposition~\ref{prop:ducktape} and partial summation.

\section{Remarks}

For an irrational number $\alpha$ in the range $0<\alpha<1$, it is
clear that the Beatty sequence $\cB_{\alpha,\beta}$ contains all
prime numbers.  In this case, since $\psi^{(\kappa)}(x)=1$ for all
$\kappa\ge 1$ and $x\in\R$, the statements in
Theorems~\ref{thm:main-two} and~\ref{thm:main-many} are consistent
with known results for the number of representations of an integer
$N$ as a sum of $\kappa$ prime numbers.

It would be interesting to see whether the results of this paper
can be extended to include irrational numbers $\alpha$ of infinite
type (with a weakened error term).

Given a sequence of real numbers $\beta_1,\ldots,\beta_\kappa$,
the techniques and results of this paper can be easily extended to
derive estimates for the number of representations of an integer
$N\equiv\kappa\pmod 2$ as a sum of $\kappa$ prime numbers,
$N=p_1+\cdots+p_\kappa$, where $p_j$ lies in the Beatty sequence
$\cB_{\alpha,\beta_j}$ for $j=1,\ldots,\kappa$.  On the other
hand, for a sequence $\alpha_1,\ldots,\alpha_\kappa$ of irrational
numbers greater than one, it appears to be much more difficult to
estimate the number of representations of an integer
$N\equiv\kappa\pmod 2$ as a sum of $\kappa$ prime numbers,
$N=p_1+\cdots+p_\kappa$, where $p_j$ lies in the Beatty sequence
$\cB_{\alpha_j,\beta_j}$ for $j=1,\ldots,\kappa$.

Finally, we have observed an interesting phenomenon.  If
$\alpha,\beta,\beta'\in\R$ with $\alpha>1$ and $\alpha$ is an
irrational number of finite type, put
$$
\cG_\kappa(\alpha,\beta;N)
=\sum_{\substack{n_1+\cdots+n_\kappa=N\\
n_1,\ldots,n_\kappa\in\cB_{\alpha,\beta}}}
\Lambda(n_1)\cdots\Lambda(n_\kappa)
$$
as before, and let $\cG_\kappa(\alpha,\beta';N)$ be defined
similarly.  If $\beta'=\beta+\alpha/\kappa$ for some fixed
$\kappa>\alpha$, then it is easy to see that the Beatty sequences
$\cB_{\alpha,\beta}$ and $\cB_{\alpha,\beta'}$ contain
\emph{different} sets of primes. Nevertheless, by
Theorem~\ref{thm:main-many} one can immediately conclude that
$$
\cG_\kappa(\alpha,\beta;N)\sim\cG_\kappa(\alpha,\beta';N)
\qquad(N\to\infty).
$$


\begin{thebibliography}{99}

\bibitem{BalPer}
A.~Balog and A.~Perelli, `Exponential sums over primes in an
arithmetic progression',  \emph{Proc.\ Amer.\ Math.\ Soc.}
\textbf{93} (1985), 578--582.

\bibitem{BaSh} W.~Banks and I.~Shparlinski,
`Prime numbers with Beatty sequences,' preprint, 2006.

\bibitem{Bug}
Y.~Bugeaud, \emph{Approximation by algebraic numbers}. Cambridge
Tracts in Mathematics, \textbf{160}. Cambridge University Press,
Cambridge, 2004.

\bibitem{Hua} L.~K.~Hua,
\emph{Introduction to Number Theory}. Springer-Verlag, Berlin
Heidelberg New York 1982.

\bibitem{Hux}
M.~N.~Huxley, \emph{The distribution of prime numbers. Large
sieves and zero-density theorems}. Clarendon Press, Oxford, 1972.

\bibitem{IwKow}
H.~Iwaniec and E.~Kowalski, \emph{Analytic number theory}. American
Mathematical Society Colloquium Publications, \textbf{53}. American
Mathematical Society, Providence, RI, 2004.

\bibitem{Khin}
A.~Y.~Khinchin, `Zur metrischen Theorie der diophantischen
Approximationen', \emph{Math.\ Z.} \textbf{24} (1926), no.~4,
706--714.

\bibitem{Koro}
N.~M.~Korobov, \emph{Exponential sums and their applications}.
Mathematics and its Applications (Soviet Series), \textbf{80}.
Kluwer Academic Publishers Group, Dordrecht, 1992.

\bibitem{KuNi}
L.~Kuipers and H.~Niederreiter, \emph{Uniform distribution of
sequences}. Pure and Applied Mathematics. Wiley-Interscience, New
York-London-Sydney, 1974.

\bibitem{Lav}
A.~F.~Lavrik, `Analytic method of estimates of trigonometric sums
by the primes of an arithmetic progression', (Russian)
\emph{Dokl.\ Akad.\ Nauk SSSR}  \textbf{248}  (1979), no.~5,
1059--1063.

\bibitem{Norton}
K.~K.~Norton, `Upper bounds for sums of powers of divisor
functions',  \emph{J.\ Number Theory}  \textbf{40}  (1992),
no.~1, 60--85.

\bibitem{Roth1}
K.~F.~Roth, `Rational approximations to algebraic numbers',
\emph{Mathematika} \textbf{2} (1955), 1--20.

\bibitem{Roth2}
K.~F.~Roth, `Corrigendum to ``Rational approximations to algebraic
numbers''',  \emph{Mathematika} \textbf{2} (1955), 168.

\bibitem{Schm}
W.~M.~Schmidt, \emph{Diophantine approximation}. Lecture Notes in
Mathematics, \textbf{785}. Springer, Berlin, 1980.

\bibitem{Vin} I.~M.~Vinogradov,
\emph{The method of trigonometrical sums in the theory of numbers}.
Dover Publications, Inc., Mineola, NY, 2004.

\bibitem{Wils}
B.~M.~Wilson, `Proofs of some formulae enunciated by Ramanujan',
\emph{Proc.\ London Math.\ Soc.} \textbf{21} (1922), 235--255.

\bibitem{Wirs} E.~Wirsing,
`Das asymptotische Verhalten von Summen \"uber multiplikative
Funktionen'. {\it Math.\ Ann.} {\bf 143} (1961), 75--102.

\end{thebibliography}
\end{document}